\documentclass[12pt,reqno]{amsart}
\usepackage{amssymb}
\usepackage{amsxtra}
\usepackage[mathscr]{eucal}

\DeclareMathOperator*{\Cat}{\displaystyle\mathbf{Cat}}

\theoremstyle{plain}
\newtheorem{Thm}{Theorem}

\newtheorem{Lem}{Lemma}
\theoremstyle{definition}
\newtheorem{Def}{Definition}
\theoremstyle{remark}

\newtheorem*{Notn}{Notation}
\newtheorem{Ex}{Example}

\numberwithin{equation}{section}

\begin{document}

\title[Partition Identities]{Partition Identities for the Multiple
Zeta Function}

\date{\today}

\author{David~M. Bradley}
\address{Department of Mathematics \& Statistics\\
         University of Maine\\
         5752 Neville Hall
         Orono, Maine 04469-5752\\
         U.S.A.}
\email[David~M. Bradley]{bradley@math.umaine.edu,
dbradley@member.ams.org}



\begin{abstract}
   We define a class of expressions for the multiple zeta function, and
   show how to determine whether an expression
   in the class vanishes identically.  The class of such identities,
   which we call partition identities, is shown to
   coincide with the class of identities
   that can be derived as a consequence of the
   stuffle multiplication rule for multiple zeta values.
\end{abstract}

\maketitle

\section{Introduction}\label{sect:Intro} For positive integer $n$
and real $s_j\ge 1$ $(j=1,2,\dots,n)$ the multiple zeta function
may be defined by
\begin{equation}
   \zeta(s_1,s_2,\dots,s_n) = \sum_{k_1>k_2>\cdots>k_n>0}\;\prod_{j=1}^n
   k_j^{-s_j}.
\label{mzvdef}
\end{equation}
The nested sum~\eqref{mzvdef} is over all positive integers
$k_1,\dots,k_n$ satisfying the indicated inequalities, and is
finite if and only if $s_1>1$ also holds.  An elementary property
of the multiple zeta function is that it satisfies the so-called
stuffle multiplication rule~\cite{BBBLa}: If
$\vec{u}=(u_1,\dots,u_m)$ and $\vec{v}=(v_1,\dots,v_n)$, then
\begin{equation}
   \zeta(\vec{u})\zeta(\vec{v})
   = \sum_{\vec{w}\in\vec{u} *\vec{v}}\zeta(\vec{w}),
\label{StuffleMult}
\end{equation}
where $\vec{u} *\vec{v}$ is the multi-set of
size~\cite{BowBradSurvey}
\[
   |\vec{u}*\vec{v}| =
   \sum_{k=0}^{\min(m,n)}\binom{m}{k}\binom{n}{k}2^k
\]
defined by the recursion
\begin{multline*}
   (s,\vec{u})*(t,\vec{v}) = \{(s,\vec{w}) : \vec{w}
   \in \vec{u}*(t,\vec{v})\} \cup \{(t,\vec{w}):
   \vec{w}\in (s,\vec{u})*\vec{v}\}\\
   \qquad\qquad\qquad\qquad\cup\{(s+t,\vec{w}):\vec{w}\in
   \vec{u}*\vec{v}\},
\end{multline*}
with initial conditions $\vec{u}*()=()*\vec{u}=\vec{u}.$ Thus, for
example,
\begin{align*}
   (s,u)*(t,v) = &\{(s,u,t,v), (s,u+t,v), (s,t,u,v), (s,t,u+v),
     (s,t,v,u)\}\\
     \cup &\{(t,s,u,v), (t,s,u+v), (t,s,v,u), (t,s+v,u),
     (t,v,s,u)\}\\
     \cup &\{(s+t,u,v), (s+t,u+v), (s+t,v,u)\},
\end{align*}
and correspondingly, we have the stuffle identity
\begin{align*}
   \zeta(s,u)\zeta(t,v) &=
   \zeta(s,u,t,v)+\zeta(s,u+t,v)+\zeta(s,t,u,v)+\zeta(s,t,u+v)+\zeta(s,t,v,u)\\
   &+\zeta(t,s,u,v)+\zeta(t,s,u+v)+\zeta(t,s,v,u)+\zeta(t,s+v,u)+\zeta(t,v,s,u)\\
   &+\zeta(s+t,u,v)+\zeta(s+t,u+v)+\zeta(s+t,v,u).
\end{align*}
The sum on the right hand side of equation~\eqref{StuffleMult}
accounts for all possible interlacings of the summation indices
when the two nested series on the left are multiplied.

In this paper, we consider a certain class of expressions (``legal
expressions'') for the multiple zeta function, consisting of a
finite linear combination of terms.  Roughly speaking, a term is a
product of multiple zeta functions, each of which is evaluated at
a sequence of sums selected from a common argument list
$(s_1,\dots,s_n)$ in such a way that each variable $s_j$ appears
exactly once in each term.  A more precise definition is given in
Section~\ref{sect:Defs}.  Once the legal expressions have been
defined, we consider the problem of determining when a legal
expression vanishes identically.  For reasons which will become
clear, we call such identities \emph{partition identities}.  It
will be seen that the problem of verifying or refuting an alleged
partition identity reduces to finite arithmetic over a polynomial
ring. Alternatively, one can first rewrite any legal expression as
a sum of single multiple zeta functions by applying the stuffle
multiplication rule to each term.  As we shall see, it is then
easy to determine whether or not the original expression vanishes
identically.

\section{Definitions}\label{sect:Defs}
Our definition of a partition identity makes use of the concept of
a set partition.  It is helpful to distinguish between set
partitions that are ordered and those that are unordered.

\begin{Def}[Unordered Set Partition]
Let $S$ be a finite non-empty set.  An \emph{unordered} set
partition of $S$ is a finite non-empty set $\mathscr{P}$ whose
elements are disjoint non-empty subsets of $S$ with union $S$.
That is, there exists a positive integer $m=|\mathscr{P}|$ and
non-empty subsets $P_1,\dots,P_m$ of $S$ such that
$\mathscr{P}=\{P_1,\dots,P_m\}$, $S=\cup_{k=1}^m P_k$, and
$P_j\cap P_k$ is empty if $j\ne k$.
\end{Def}

\begin{Def}[Ordered Set Partition]
Let $S$ be a finite non-empty set.  An \emph{ordered} set
partition of $S$ is a finite ordered tuple $\vec P$ of disjoint
non-empty subsets of $S$ such that the union of the components of
$\vec P$ is equal to $S$. That is, there exists a positive integer
$m$ and non-empty subsets $P_1,\dots,P_m$ of $S$ such that $\vec
P$ can be identified with the ordered $m$-tuple $(P_1,\dots,P_m)$,
$\cup_{k=1}^m P_k=S$, and $P_j\cap P_k$ is empty if $j\ne k$.
\end{Def}

\begin{Def}[Legal Term]\label{def:legalterm}
Let $n$ be a positive integer and let $\vec{s} = (s_1,\dots,s_n)$
be an ordered tuple of $n$ real variables with $s_j>1$ for $1\le
j\le n$. Let ${\mathscr{P}}=\{P_1,\dots,P_m\}$ be an unordered set
partition of the first $n$ positive integers $\{1,2,\dots,n\}$.
For each positive integer $k$ such that $1\le k\le m$, let
$\vec{P_k}=(P_k^{(1)},P_k^{(2)},\dots,P_k^{(\alpha_k)})$ be an
ordered set partition of $P_k$, and let
\[
   t_k^{(r)} = \sum_{j\in P_k^{(r)}} s_j,
   \qquad 1\le r\le \alpha_k=|\vec P_k|.
\]
A \emph{legal term} for $\vec{s}$ is a product of the form
\begin{equation}
   \prod_{k=1}^m
   \zeta\big(t_k^{(1)},t_k^{(2)},\dots,t_k^{(\alpha_k)}\big),
\label{LegalTerm}
\end{equation}
and every legal term for $\vec{s}$ has the form~(\ref{LegalTerm})
for some unordered set partition ${\mathscr{P}}$ of
$\{1,2,\dots,n\}$ and ordered subpartitions $\vec{P_k}$, $1\le
k\le |{\mathscr{P}}|$.
\end{Def}

\begin{Ex} The product
$\zeta(s_6,s_2+s_5,s_1+s_8+s_9)\zeta(s_3+s_4,s_{10})\zeta(s_7)$ is
a legal term for the 10-tuple
$(s_1,s_2,s_3,s_4,s_5,s_6,s_7,s_8,s_9,s_{10})$ arising from the
partition $\{P_1,P_2,P_3\}$ of the set $\{1,2,3,4,5,6,7,8,9,10\}$,
where $P_1=\{1,2,5,6,8,9\}$ has ordered subpartition $\vec
P_1=(\{6\},\{2,5\},\{1,8,9\})$, $P_2=\{3,4,10\}$ has ordered
subpartition $\vec P_2=(\{3,4\},\{10\})$, and $P_3=\{7\}$ has
ordered subpartition $\vec P_3=(\{7\}).$
\end{Ex}

\begin{Def}[Legal Expression]
Let $n$ be a positive integer, and let $\vec{s}=(s_1,\dots,s_n)$
be an ordered tuple of $n$ real variables with $s_j>1$ for $1\le
j\le n$. A \emph{legal expression} for $\vec{s}$ is a finite
$\mathbf{Z}$-linear combination of legal terms for $\vec{s}$. That
is, for any positive integer $q$, integers $a_h$, and legal terms
$T_h$ for $\vec{s}$ $(1\le h\le q)$, the sum $\sum_{h=1}^q a_h
T_h$ is a legal expression for $\vec{s}$, and every legal
expression for $\vec{s}$ has this form.
\end{Def}

\begin{Def}[Partition Identity]
A \emph{partition identity} is an equation of the form
$\mathrm{LHS}=0$ for which there exists a positive integer $n$ and
real variables $s_j>1$ $(j=1,2,\dots,n)$ such that $\mathrm{LHS}$
is a legal expression for $(s_1,\dots,s_n)$, and the equation
holds true for all real values of the variables $s_j>1$.
\end{Def}

\begin{Ex}\label{Ex:Partn} The equation
\begin{multline*}
   2\zeta(s_1+s_2+s_3)-\zeta(s_2)\zeta(s_1+s_3)-\zeta(s_3)\zeta(s_1+s_2)+\zeta(s_1+s_2,s_3)\\
   +\zeta(s_2,s_1+s_3)+\zeta(s_1+s_3,s_2)+\zeta(s_3,s_1+s_2)=0
\end{multline*}
is a partition identity, and is easily verified by expanding the
two products $\zeta(s_2)\zeta(s_1+s_3)$ and
$\zeta(s_3)\zeta(s_1+s_2)$ using the stuffle multiplication
rule~\eqref{StuffleMult} and then collecting multiple zeta
functions with identical arguments. A natural question is whether
\emph{every} partition identity can be verified in this way.  We
provide an affirmative answer to this question in
Section~\ref{sect:canonical}.  An alternative method for verifying
partition identities is given in Section~\ref{sect:rational}.
\end{Ex}

\section{Rational Functions}\label{sect:rational}
Here, we describe a method by which one can determine whether or
not a legal expression vanishes identically, or equivalently,
whether or not an alleged partition identity is in fact a true
identity. It will be seen that the problem reduces to that of
checking whether or not an associated rational function identity
is true.  This latter check can be accomplished in a completely
deterministic and mechanical fashion by clearing denominators and
expanding the resulting multivariate polynomials. More
specifically, we associate rational functions with legal terms in
such a way that the alleged partition identity holds if and only
if the corresponding rational function identity, in which each
legal term is replaced by its associated rational function, holds.
The rational function corresponding to~(\ref{LegalTerm}) is the
function of $n$ real variables $x_1>1,\dots,x_n>1$ defined by
\begin{equation}
   R(x_1,x_2,\dots,x_n):=\prod_{k=1}^m \prod_{\beta=1}^{\alpha_k}
   \bigg(\prod_{\lambda=1}^{\beta}\prod_{j\in P_k^{(\lambda)}}
   x_j-1\bigg)^{-1}.
\label{RatTerm}
\end{equation}

\begin{Thm}\label{Thm:Rat} Let $q$ be a positive integer, and let
$E=\sum_{h=1}^q a_h T_h$ be a legal expression for $\vec
s=(s_1,\dots,s_n)$ \textup{(}i.e.\ each $a_h\in{\mathbf Z}$ and
$T_h$ is a legal term for $\vec s$, $1\le h\le q$\textup{)}. Let
$L=\sum_{h=1}^q a_h r_h$ be the expression obtained by replacing
each legal term $T_h$ by its corresponding rational function
according to the rule that associates~\eqref{RatTerm}
with~\eqref{LegalTerm}. Then $E$ vanishes identically if and only
if $L$ does.
\end{Thm}

\begin{Ex}
The rational function identity which Theorem~\ref{Thm:Rat} asserts
is equivalent to the partition identity of Example~\ref{Ex:Partn}
is
\begin{multline*}
   \frac{2}{x_1x_2x_3-1}-\frac{1}{x_2-1}\cdot\frac{1}{x_1x_3-1}
   -\frac{1}{x_3-1}\cdot\frac{1}{x_1x_2-1}+\frac{1}{(x_1x_2-1)(x_1x_2x_3-1)}\\
   +\frac{1}{(x_2-1)(x_1x_2x_3-1)}+\frac{1}{(x_1x_3-1)(x_1x_2x_3-1)}
   +\frac{1}{(x_3-1)(x_1x_2x_3-1)}=0,
\end{multline*}
which can be readily verified by hand, or with the aid of a
suitable computer algebra system.
\end{Ex}

\noindent{\bf Proof of Theorem~\ref{Thm:Rat}.} It is immediate
from the partition integral~\cite{BBBLa} representation for the
multiple zeta function that every legal term on $(s_1,\dots,s_n)$
is an $n$-dimensional integral transform of its associated
rational function multiplied by the common kernel $\prod_{j=1}^n
(\log x_j)^{s_j-1}/\Gamma(s_j)x_j$. Explicitly,
\begin{multline*}
   \prod_{k=1}^m
   \zeta\bigg(\sum_{j\in P_k^{(1)}} s_j, \sum_{j\in P_k^{(2)}} s_j,
   \dots,\sum_{j\in P_k^{(\alpha_k)}}
   s_j\bigg)\\
   = \int_1^\infty\cdots\int_1^\infty \bigg\{\prod_{k=1}^m
   \prod_{\beta=1}^{\alpha_k}\bigg(\prod_{\lambda=1}^\beta
   \prod_{j\in P_k^{(\lambda)}} x_j -1 \bigg)^{-1}\bigg\}
   \prod_{j=1}^n \frac{(\log x_j)^{s_j-1}}{\Gamma(s_j)x_j}\,dx_j.
\end{multline*}
Linearity of the integral implies that if $L\equiv 0$ then
$E\equiv 0$. The real content of Theorem~\ref{Thm:Rat} is that
converse the also holds.  To prove this, we first note that the
rational function~\eqref{RatTerm} is continuous on the $n$-fold
Cartesian product of open intervals
$(1..\infty)^n=\{(x_1,\dots,x_n)\in{\mathbf R}^n : x_j>1, 1\le
j\le n\}$ and $|R(x_1,\dots,x_n)\prod_{j=1}^n x_j|$ is bounded on
any $n$-fold Cartesian product of half-open intervals of the form
$[c..\infty)^n=\{(x_1,\dots,x_n)\in{\mathbf R}^n : x_j\ge c, 1\le
j\le n\}$ with $c>1$. These properties obviously extend to linear
combinations of rational functions of the form~\eqref{RatTerm},
and thus to complete the proof of Theorem~\ref{Thm:Rat}, it
suffices to establish the following result.

\begin{Lem} Let $n$ be a positive integer and let $R$ be a
continuous real-valued function of $n$ real variables defined on
the $n$-fold Cartesian product of open intervals $(1..\infty)^n$.
Suppose there exists a constant $c>1$ such that
$|R(x_1,x_2,\ldots,x_n)\prod_{j=1}^n x_j|$ is bounded on the
$n$-fold Cartesian product of half-open intervals $[c..\infty)^n$.
Suppose further that there exist non-negative real numbers
$s_1^*,s_2^*,\dots,s_n^*$ such that the $n$-dimensional multiple
integral
\[
   \int_1^\infty\cdots\int_1^\infty  R(x_1,x_2,\ldots,x_n)
   \prod_{j=1}^n \left(\log x_j\right)^{s_j}
   \frac{dx_j}{x_j}
\]
vanishes whenever $s_j>s_j^*$ for $1\le j\le n.$  Then $R$
vanishes identically.
\end{Lem}

\noindent{\bf Proof.} Fix $s_j>s_j^*$ for $1\le j\le n$.  Let
$T:[1..\infty)\to {\mathbf R}$ be given by the convergent
$(n-1)$-dimensional multiple integral
\[
   T(x) := \int_1^\infty\cdots\int_1^\infty R(x_1,\ldots,x_{n-1},x)
   \prod_{j=1}^{n-1}\left(\log x_j\right)^{s_j}
   \frac{dx_j}{x_j} .
\]
Then $T(x)=O(1/x)$ as $x\to\infty$.  It follows that the Laplace
Transform $$
   F(z) := \int_0^\infty e^{-zu}\,T(e^u)\,du
$$ is analytic in the right half-plane $\{z\in{\mathbf C}:
\Re(z)>-1\}$, and for all positive integers $m>s_n^*$,
\[
   \left(-\frac{d}{dz}\right)^mF(z)\bigg|_{z=0}
   =\int_0^\infty u^m\,T(e^u)\,du
   =\int_1^\infty\left(\log x\right)^m T(x)\,\frac{dx}{x}
   = 0.
\]
By Taylor's theorem, $F$ is a polynomial.  Letting $z\to+\infty$
in the definition of $F$, we see that in fact, $F$ must be the
zero polynomial.  By the uniqueness theorem for Laplace transforms
(see eg.~\cite{Widd}), the set of $x>1$ for which $T(x)\ne 0$ is
of Lebesgue measure zero.  Since $T$ is continuous, it follows
that $T(x)=0$ for all $x>1$.  If $n=1$, then $T=R$ and we're done.
Otherwise, fix $x>1$, and suppose the result holds for $n-1$.
Since in the above argument, $s_1>s_1^*,\dots,s_{n-1}>s_{n-1}^*$
were arbitrary, $T(x)=0$ implies $R(x_1,x_2,\dots,x_{n-1},x)=0$
for all $x_1,x_2,\dots,x_{n-1}$ by the inductive hypothesis. Since
this is true for each fixed $x>1$, the result follows. \qed

\section{Stuffles and Partition Identities}\label{sect:canonical}

As in~\cite{BBBLa}, we define the class of stuffle identities to
be the set of all identities of the form~\eqref{StuffleMult}.
In~\cite{BBBLa}, it is shown that every stuffle identity is a
consequence of a corresponding rational function identity.  In the
previous section of the present paper, using a different method of
proof, we established the more general result that every partition
identity is a consequence of a corresponding rational function
identity. Clearly every stuffle identity is a partition identity,
but not conversely. Nevertheless, we shall see that every
partition identity is a consequence of the stuffle multiplication
rule. More specifically, we provide an affirmative answer to the
question raised at the end of Example~\ref{Ex:Partn} in
Section~\ref{sect:Defs}.

\begin{Notn}
We introduce the concatenation operator $\Cat$, which will be
useful for expressing argument sequences without recourse to
ellipses. For example, $\Cat_{k=1}^m t_j$ denotes the sequence
$t_1,\dots,t_m$.
\end{Notn}

As we noted previously, by applying the stuffle multiplication
rule~\eqref{StuffleMult} to legal terms, any legal expression on
$(s_1,\dots,s_n)$ can be rewritten as a finite ${\mathbf
Z}$-linear combination of single multiple zeta functions of the
form
\[
   \sum_{h=1}^q a_h\, \zeta\big(\Cat_{k=1}^{\alpha_h} \sum_{j\in
   P_k}s_j\big)=
   \sum_{h=1}^q a_h\, \zeta\big(\sum_{j\in P_1}s_j,
   \sum_{j\in P_2}s_j,\dots,\sum_{j\in P_{\alpha_h}}
   s_j\big),
\]
where the coefficients $a_h$ are integers, $q$ is a positive
integer, and $(P_1,\dots,P_{\alpha_h})$ is an ordered set
partition of the first $n$ positive integers $\{1,2,\dots,n\}$ for
each $h=1,2,\dots,q$. Thus, it suffices to prove the following
result.

\begin{Thm}\label{Thm:Canonical} Let $F$ be a finite non-empty
set of positive integers, and let $\{s_j : j\in F\}$ be a set of
real variables, each exceeding 1.  Suppose that for all $s_j>1$,
\begin{equation}
   \sum_{\vec P\vDash F} c_{\vec P}\, \zeta\big( \Cat_{k=1}^{|\vec P|}
   \sum_{j\in P_k} s_j\big) = 0,
\label{independent}
\end{equation}
where the sum is over all ordered set partitions $\vec P$ of $F$,
$P_k$ denotes the $k^{\mathrm{th}}$ component of $\vec P$, and the
coefficients $c_{\vec P}$ are real numbers depending only on $\vec
P$.  Then each $c_{\vec P}=0$.
\end{Thm}

\noindent{\bf Proof.}  We argue by induction on the cardinality of
$F$, the case $|F|=1$ being trivial.  To clarify the argument, we
present the cases $|F|=2$ and $|F|=3$ before proceeding to the
inductive step.

When $|F|=2$, the identity ~\eqref{independent} takes the form
\[
   c_1\zeta(s,t)+c_2\zeta(t,s)+c_3\zeta(s+t)=0,\qquad s>1,\; t>1.
\]
By Theorem~\ref{Thm:Rat}, this is equivalent to the rational
function identity
\[
   \frac{c_1}{(x-1)(xy-1)}+\frac{c_2}{(y-1)(xy-1)}+\frac{c_3}{xy-1}=0,
   \qquad x>1,\; y>1.
\]
Letting $x\to1+$ shows that we must have
\[
   \frac{c_1}{y-1} = 0 \quad\Longrightarrow\quad c_1\zeta(t)=0
   \quad\Longrightarrow \quad c_1=0.
\]
Similarly, letting $y\to1+$ shows that $c_2=0$.  Since
the remaining term must vanish, we must have $c_3=0$ as well.

When $|F|=3$, the identity~\eqref{independent} takes the form
\begin{align*}
  0&=c_1\zeta(s,t,u)+c_2\zeta(s,u,t)+c_3\zeta(t,s,u)+c_4\zeta(t,u,s)+c_5\zeta(u,s,t)
   +c_6\zeta(u,t,s)\\
   &+c_7\zeta(s,t+u)+c_8\zeta(t+u,s)+c_9\zeta(t,s+u)+c_{10}\zeta(s+u,t)
   +c_{11}\zeta(u,s+t)\\
   &+c_{12}\zeta(s+t,u)+c_{13}\zeta(s+t+u), \qquad s>1,\; t>1,\; u>1,
\end{align*}
which, by Theorem~\ref{Thm:Rat}, is equivalent to the rational
function identity
\begin{align}
  0 &= \frac{c_1}{(x-1)(xy-1)(xyz-1)}+\frac{c_2}{(x-1)(xz-1)(xyz-1)}
   +\frac{c_3}{(y-1)(xy-1)(xyz-1)}\nonumber\\
   &+\frac{c_4}{(y-1)(yz-1)(xyz-1)}
   +\frac{c_5}{(z-1)(xz-1)(xyz-1)}+\frac{c_6}{(z-1)(yz-1)(xyz-1)}\nonumber\\
   &+\frac{c_7}{(x-1)(xyz-1)}+\frac{c_8}{(yz-1)(xyz-1)}
   +\frac{c_9}{(y-1)(xyz-1)}+\frac{c_{10}}{(xz-1)(xyz-1)}\nonumber\\
   &+\frac{c_{11}}{(z-1)(xyz-1)}+\frac{c_{12}}{(xy-1)(xyz-1)}
   +\frac{c_{13}}{xyz-1},\qquad x>1,\; y>1,\; z>1.
   \label{canon3}
\end{align}
If we let $x\to1+$ in~\eqref{canon3}, then for the singularities
to cancel, we must have
\[
   0=\frac{c_1}{(y-1)(yz-1)}+\frac{c_2}{(z-1)(yz-1)}+\frac{c_7}{yz-1},
   \qquad y>1,\; z>1,
\]
which, in light of Theorem~\ref{Thm:Rat}, implies the identity
\[
   0=c_1\zeta(t,u)+c_2\zeta(u,t)+c_7\zeta(t+u), \qquad t>1,\; u>1.
\]
Having proved the $|F|=2$ case of our result, we see that this
implies $c_1=c_2=c_7=0$.  Similarly, letting $y\to 1+$
in~\eqref{canon3} gives $c_3=c_4=c_9=0$, and letting $z\to 1+$
in~\eqref{canon3} gives $c_5=c_6=c_{11}=0$.  At this point, only 4
terms in~\eqref{canon3} remain.  Letting $yz\to1+$ now shows that
$c_8=0$, and the remaining coefficients can be shown to vanish
similarly.

For the inductive step, let $|F|>1$ and suppose
Theorem~\ref{Thm:Canonical} is true for all non-empty sets of
positive integers of cardinality less than $|F|$.  Suppose also
that~\eqref{independent} holds. By Theorem~\ref{Thm:Rat}, it
follows that the rational function identity
\begin{equation}
   \sum_{\vec P\vDash F} c_{\vec P}\prod_{m=1}^{|\vec P|}
   \big(\prod_{k=1}^m \prod_{j\in P_k} x_j-1\big)^{-1} =0
\label{canon}
\end{equation}
holds for all $x_j>1$, $j\in F$.  Fix $f\in F$ and let $x_f\to 1+$
in~\eqref{canon}. For the singularities to cancel, we must have
\[
   \sum_{\substack{\vec P\vDash F\\ P_1=\{f\}}}
   c_{\vec P} \prod_{m=2}^{|\vec P|}\big(\prod_{k=2}^m \prod_{j\in
   P_k} x_j-1\big)^{-1}=0,
\]
which, by Theorem~\ref{Thm:Rat}, implies that
\[
    \sum_{\substack{\vec P\vDash F\\ P_1=\{f\}}}
    c_{\vec P}\, \zeta\big(\Cat_{m=2}^{|\vec P|} \sum_{j\in
    P_m}s_j\big)
    =
    \sum_{\vec P\vDash F\setminus\{f\}}
    c_{\vec P}\, \zeta\big(\Cat_{m=1}^{|\vec P|} \sum_{j\in
    P_m}s_j\big)
    =0.
\]
By the inductive hypothesis, $c_{\vec P}=0$ for every ordered set
partition $\vec P$ of $F$ whose first component $P_1$ is the
singleton $\{f\}$.  Since $f\in F$ was arbitrary, it follows that
$c_{\vec P}=0$ for every ordered set partition $\vec P$ of $F$
whose first component $P_1$ consists of a single element.

Proceeding inductively, suppose we've shown that $c_{\vec P}=0$
for every ordered set partition $\vec P$ of $F$ with
$|P_1|=r-1<|F|$.  Let $G$ be a subset of $F$ of cardinality $r$.
If in~\eqref{canon} we now let $x_g\to1+$ for each $g\in G$, then
as the singularities in the remaining terms~\eqref{canon} must
cancel, we must have
\[
   \sum_{\substack{ \vec P\vDash F\\ P_1=G}}
   c_{\vec P}\prod_{m=2}^{|\vec P|} \big(\prod_{k=2}^m
   \prod_{j\in P_k}x_j-1\big)^{-1}=0.
\]
Theorem~\ref{Thm:Rat} then implies that
\[
   \sum_{\substack{ \vec P\vDash F\\ P_1=G}}
   c_{\vec P}\,\zeta\big(\Cat_{m=2}^{|\vec P|}\sum_{j\in P_m}s_j\big)
   =\sum_{\vec P\vDash F\setminus G}
   c_{\vec P}\,\zeta\big(\Cat_{m=1}^{|\vec P|}\sum_{j\in P_m}s_j\big)
   = 0.
\]
By the inductive hypothesis, $c_{\vec P}=0$ for every ordered set
partition $\vec P$ of $F$ with first component equal to $G$. Since
$G$ was an arbitrary subset of $F$ of cardinality $r$, it follows
that $c_{\vec P}=0$ for every ordered set partition $\vec P$ of
$F$ whose first component has cardinality $r$.  By induction on
$r$, it follows that $c_{\vec P}=0$ for every ordered set
partition $\vec P$ of $F$, as claimed. \qed

Since every partition identity is thus a consequence of the
stuffle multiplication rule, it follows that any multi-variate
function that obeys a stuffle multiplication rule will satisfy
every partition identity satisfied by the multiple zeta function.
For example, suppose we fix a positive integer $N$ and a set $F$
of functions $f:{\mathbf Z}^{+}\to {\mathbf C}$ closed under
point-wise addition. For positive integer $n$ and
$f_1,\dots,f_n\in F$, define
\[
   z_N(f_1,f_2,\dots,f_n) := \sum_{N>k_1>k_2>\cdots>k_n>0}\;
   \prod_{j=1}^n \exp(f_j(k_j)),
\]
where the sum is over all positive integers $k_1,\dots,k_n$
satisfying the indicated inequalities.  Then for all $g,h\in F$,
we have $z_N(g)z_N(h)=z_N(g,h)+z_N(h,g)+z_N(g+h)$, and more
generally, if $\vec{g}=(g_1,\dots,g_m)$ and
$\vec{h}=(h_1,\dots,h_n)$ are vectors of functions in $F$, then
$z_N$ obeys the stuffle multiplication rule
\[
   z_N(\vec{g})z_N(\vec{h})
   = \sum_{\vec{f}\in\vec{g} *\vec{h}}z_N(\vec{f}).
\]
We assert that $z_N$ satisfies every partition identity satisfied
by $\zeta$. For example, let $n$ be a positive integer, and let
${\mathfrak{S}_n}$ denote the group of $n!$ permutations of the
first $n$ positive integers $\{1,2,\dots,n\}$. Let $s_1,\dots,s_n$
be real variables, each exceeding 1.  Using a counting argument,
Hoffman~\cite{Hoff1} proved the partition identity
\begin{equation}
   \sum_{\sigma\in{\mathfrak{S}_n}}
   \zeta\big(\Cat_{k=1}^n s_{\sigma(k)}\big)
   = \sum_{{\mathscr{P}} \vdash \{1,\dots,n\}}
   (-1)^{n-|{\mathscr{P}}|}\prod_{P\in{\mathscr{P}}}
   (|P|-1)!\,\zeta\big(\sum_{j\in P}s_j\big),
\label{hoff}
\end{equation}
in which the sum on the right extends over all unordered set
partitions ${\mathscr{P}}$ of the first $n$ positive integers
$\{1,2,\dots,n\}$, and of course $|{\mathscr{P}}|$ denotes the
number of parts in the partition ${\mathscr{P}}$.  In light of
Theorem~\ref{Thm:Canonical}, it follows that Hoffman's
identity~\eqref{hoff} depends on only the stuffle multiplication
property~\eqref{StuffleMult} of the multiple zeta function; whence
any function satisfying a stuffle multiplication rule will also
satisfy~\eqref{hoff}.  In particular, with $z_N$ defined as above,
\[
   \sum_{\sigma\in{\mathfrak{S}_n}}
   z_N\big(\Cat_{k=1}^n f_{\sigma(k)}\big)
   = \sum_{{\mathscr{P}} \vdash \{1,\dots,n\}}
   (-1)^{n-|{\mathscr{P}}|}\prod_{P\in{\mathscr{P}}}
   (|P|-1)!\,z_N\big(\sum_{j\in P}f_j\big).
\]
Finally, we note that by Theorem~\ref{Thm:Rat}, the rational
function identity
\[
   \sum_{\sigma\in{\mathfrak{S}_n}} \prod_{k=1}^n
   \big(\prod_{j=1}^k x_{\sigma(j)} -1\big)^{-1}
   =\sum_{{\mathscr{P}}\vdash \{1,\dots,n\}}
   (-1)^{n-|{\mathscr{P}}|}\prod_{P\in{\mathscr{P}}}
   (|P|-1)!\,\big(\prod_{j\in P}x_j-1\big)^{-1}
\]
is equivalent to~\eqref{hoff}.

\end{document}